\documentclass[14]{article}
\usepackage[margin=1in]{geometry}
\usepackage{amsfonts,amsmath,amssymb}
\numberwithin{equation}{section}

\title{The Variational Iteration Method for Solving Asymmetric System of PDEs}
\author{Abdulhameed Qahtan Abbood Altai \\ 
\small University of Babylon, Babil, Iraq 51002 \\
\small ahbabil1983@gmail.com}

\begin{document}

\maketitle

\paragraph{Abstract.} We propose a method to obtain iterative schemes guarantee unique solutions for systems of partial differential equations that are not symmetric with respect to the time by generalizing He's variational iteration method and using Banach's fixed point theorem. Then, iterative schemes for systems of incompressible fluid flow and incompressible micropolar fluid flow will be created by applying the generalized He's variational iteration method.

\paragraph{Keywords:} He's variational iteration method, Navier-Stokes equations, microrotational velocity equations.

\section{Introduction} 
Variational iteration method (VIM) was proposed and used to solve many kinds of PDEs by He $[7-16]$. Its advantage is to provide a faster successive approximation of the (exact or numerical) solution comparing with Adomian's decomposition method [23,1,2] and it does not depend on small parameters comparing with perturbation methods. Moreover, it was used to solve the systems of PDEs which are symmetric with respect to the time [24]. In this article, we try to study asymmetric systems of PDEs with respect to the time; in other words, systems of PDEs for unknowns $\textit{\textbf{u}}_{k}, k = 1, \cdots, m $ and unknowns $p_{k}, k = 1, \cdots, m $, such that every PDE$_{k}$ has the Laplacian operator $\triangle \textit{\textbf{u}}_{k}$ and the operator $\nabla p_{k}$ for $k = 1, \cdots, m $, and they are symmetric for $\textit{\textbf{u}}_{k}, k = 1, \cdots, m $ and not symmetric for $p_{k}, k = 1, \cdots, m $ with respect to the time.  Our way is to reform the system of PDEs and normalize every Laplacian operator $\triangle \textit{\textbf{u}}_{k}$ in every PDE$_{k}$ for $k = 1, \cdots, m $ to make the system symmetric with respect to the location $\textit{\textbf{x}} = (x_{1}, \cdots, x_{d})$ and then we multiply every PDE$_{k}$ by the operator $\nabla$ and the unit outer normal $\textit{\textbf{n}}$ to $\partial\Omega$ to create a second system such that every PDE$_{k}$ has the Laplcian operator $\triangle p_{k}$ for unknowns $p_{k}, k = 1, \cdots, m$, and it is symmetric with respect to the location $\textit{\textbf{x}} = (x_{1}, \cdots, x_{d})$ as well. Hence, we can apply the VIM to both systems together in any  $x_{l}-$direction we like, for $l = 1, \cdots, d$, to get iterative schemes guarantee unique solutions $ \textit{\textbf{u}}_{k}, p_{k}, k = 1, \cdots, m $ using Banach's fixed point theorem. This method can be applied for systems of incompressible fluid flow (Navier-Stokes equations) and incompersible micropolar fluid flow (Navier-Stokes equations and microrotational velocity equations). So, the outline of the paper is as follows: in section 2, we generalize He's variational iteration method to the asymmetric systems of PDEs with respect to the time. Then, we apply our technique to the systems of incompressible fluid flow in section 3, and incompersibble micropolar fluid flow in section 4.

\section{Generalized He's variational iteration method} Let $\textit{\textbf{u}}_{k} = (u_{k,1}, \cdots, u_{k,d}):\Omega (\subset \mathbb{R}^d)  \times [0,T] \to \mathbb{R}^d$ and $p_{k}:\Omega (\subset \mathbb{R}^d)  \times [0,T] \to \mathbb{R}, k = 1, \cdots, m$  and $d \in \mathbb{N}$, and consider the system of partial differential equations that is not symmetric in $p_{k}, k = 1, \cdots, m$, with respect to $t \in \left[ 0, T \right] $, of the form
\begin{align}
\begin{array}{llllll}
L_{t} {\textit{\textbf{u}}}_{1} &+  R_{1} \left( \triangle{\textit{\textbf{u}}_{1}}, \textit{\textbf{u}}_{2}, \cdots, \textit{\textbf{u}}_{m}, \nabla p_{1} \right) &+  N_{1} \left( \textit{\textbf{u}}_{1}, \cdots, \textit{\textbf{u}}_{m} \right) &= \textit{\textbf{f}}_{1} \\
& \vdots       & \vdots   & \vdots      \\
L_{t} \textit{\textbf{u}}_{m} &+  R_{m} \left( \textit{\textbf{u}}_{1}, \cdots, \textit{\textbf{u}}_{m-1}, \triangle{\textit{\textbf{u}}_{m}}, \nabla p_{m} \right) &+  N_{m} \left( \textit{\textbf{u}}_{1}, \cdots, \textit{\textbf{u}}_{m} \right) &= \textit{\textbf{f}}_{m}
\end{array} 
\end{align} 
with initial data $\textit{\textbf{u}}_{k}(\textit{\textbf{x}},0) = \textit{\textbf{u}}_{k,0}(\textit{\textbf{x}}), k = 1, \cdots, m$, where $L_{t}$ is a first-order partial differential operator, $R_{k}$ and $N_{k}, k = 1, \cdots, m$, are linear and nonlinear operators respectively, and $\textit{\textbf{f}}_{k} = (f_{k,1}, \cdots,f_{k,d}), k = 1,\cdots, m$, are source terms. Since the Laplacian terms $\triangle{\textit{\textbf{u}}_{k}}$ are in linear operators $R_{k}, k = 1, \cdots, m$, we can normalize thier coefficients and reform the system $(2.1)$ to 
\begin{align}
\begin{array}{llllll} 
\triangle{\textit{\textbf{u}}_{1}} &+  R_{1} \left( \textit{\textbf{u}}_{1}, \cdots, \textit{\textbf{u}}_{m}, \nabla p_{1} \right) &+  N_{1} \left( \textit{\textbf{u}}_{1}, \cdots, \textit{\textbf{u}}_{m} \right) &= \textit{\textbf{h}}_{1} \\
& \vdots       & \vdots   & \vdots      \\
\triangle{\textit{\textbf{u}}_{m}} &+  R_{m} \left( \textit{\textbf{u}}_{1}, \cdots, \textit{\textbf{u}}_{m} , \nabla p_{m} \right) &+  N_{m} \left( \textit{\textbf{u}}_{1}, \cdots, \textit{\textbf{u}}_{m} \right) &= \textit{\textbf{h}}_{m}
\end{array} 
\end{align} 
To create a system of partial differential equations symmetric in $p_{k}; k = 1, \cdots, m$, with respect to $\mathbf{x}$, we take the operator $\nabla$ and the unit outer normal $\textit{\textbf{n}}$ to $\partial\Omega$ on both sides of every equations in $(2.1)$, to get
\begin{align}
\begin{array}{llllll}
\triangle p_{1} &+  \nabla R_{1} \left( \textit{\textbf{u}}_{1}, \cdots, \textit{\textbf{u}}_{m} \right) &+ \nabla N_{1} \left( \textit{\textbf{u}}_{1}, \cdots, \textit{\textbf{u}}_{m} \right) &= \nabla \textit{\textbf{f}}_{1} \\
&\vdots    & \vdots     &\vdots      \\
\triangle p_{m} &+  \nabla R_{m} \left( \textit{\textbf{u}}_{1}, \cdots, \textit{\textbf{u}}_{m} \right) &+  \nabla N_{m} \left( \textit{\textbf{u}}_{1}, \cdots, \textit{\textbf{u}}_{m} \right) &= \nabla \textit{\textbf{f}}_{m}.
\end{array}  
\end{align} 
with conditions $\dfrac{\partial p_{k}}{\partial \textit{\textbf{n}}} = \left( L_{\mathbf{x}} \textit{\textbf{u}}_{k} + R_{k} + N_{k} - \textit{\textbf{f}}_{k} \right). \textit{\textbf{n}}; k = 1, \cdots, m$.
Then, for $k = 1, \cdots, m$, the correctional functionals for the system $(2.2)$ and for the system $(2.3)$ in any $x_{l}-$direction, $l = 1, \cdots, d$, are
\begin{equation}
\begin{aligned}
u_{k,1,n+1}(\mathbf{x},t) &= u_{k,1,n}(\mathbf{x},t) + \int_{\Omega_{x_{l}}} \lambda_{k,1}(x_{l}) [ \triangle{u_{k,1,n}}(\mathbf{x},\tau) +  R_{k,1} \left( \tilde{u}_{1,1,n}, \cdots, \tilde{u}_{1,d,n}, \cdots, \tilde{u}_{m,1,n}, \cdots, \tilde{u}_{m,d,n}, \nabla \tilde{p}_{k,n} \right) \\
&+  N_{k,1} \left( \tilde{u}_{1,1,n}, \cdots, \tilde{u}_{1,d,n}, \cdots, \tilde{u}_{m,1,n}, \cdots, \tilde{u}_{m,d,n} \right) - f_{k,1}(\mathbf{x},\tau) ] d x_{l} \\
& \vdots     \\
u_{k,d,n+1}(\mathbf{x},t) &= u_{k,d,n}(\mathbf{x},t) + \int_{\Omega_{x_{l}}} \lambda_{k,d}(x_{l}) [ \triangle{u_{k,d,n}}(\mathbf{x},\tau) +  R_{k,d} \left( \tilde{u}_{1,1,n}, \cdots, \tilde{u}_{1,d,n}, \cdots, \tilde{u}_{m,1,n}, \cdots, \tilde{u}_{m,d,n}, \nabla \tilde{p}_{k,n} \right) \\
&+  N_{k,d} \left( \tilde{u}_{1,1,n}, \cdots, \tilde{u}_{1,d,n}, \cdots, \tilde{u}_{m,1,n}, \cdots, \tilde{u}_{m,d,n} \right) - f_{k,d}(\mathbf{x},\tau) ] d x_{l} \\
p_{k,n+1}(\mathbf{x},t) &= p_{k,n}(\mathbf{x},t) + \int _{\Omega_{x_{l}}} \mu_{k}(x_{l}) [ \triangle p_{k,n}(\mathbf{x},t) +  \nabla R_{k} \left( \tilde{\textit{\textbf{u}}}_{1}, \cdots, \tilde{\textit{\textbf{u}}}_{k} \right) + \nabla N_{k} \left( \tilde{\textit{\textbf{u}}}_{1}, \cdots, \tilde{\textit{\textbf{u}}}_{k} \right) - \nabla \textit{\textbf{f}}_{k}(\mathbf{x},t)] d x_{l} 
\end{aligned}
\end{equation}
where $\lambda_{k,l}, \mu_{k}, k = 1, \cdots, m; l = 1, \cdots, d$, are general Lagrange multipliers, and $\tilde{u}_{k,l,n}, \tilde{p}_{k,n}, k = 1, \cdots, m; l = 1, \cdots, d$ are restricted variations, that is $\delta \tilde{u}_{k,l,n} = 0, k = 1, \cdots, m; l = 1, \cdots, d$. By the variational theory and via integration by parts, the Lagrange multipliers $\lambda_{k,l}, \mu_{k}, k = 1, \cdots, m; l = 1, \cdots, d$ can be identified. To obtain unique solution of above iterative scheme, we consider the operators $ U_{k,l,n}[u_{k,l,n}] = u_{k,l,n+1}, P_{k,n}[p_{k,n}] = p_{k,n+1}, n \geq 0, k = 1, \cdots, m; l = 1, \cdots, d$ and use Banach's fixed point theorem which state: if $X$ is a Banach space and $A:X \to X$ is a nonlinear mapping satisfying 
\begin{align*}
\Vert A[f_{1}] - A[f_{2}] \Vert \leq \gamma \Vert f_{1} - f_{2} \Vert, \ f_{1},f_{1} \in X
\end{align*}
for $\gamma \in [0,1)$, then $A$ has a unique fixed point [5]. So, the sufficient condition to approximate the iteravtive scheme obtained by the generalized He's variational iteration method is strictly contraction of $ U_{k,l,n}, P_{k,n}, n \geq 0, k = 1, \cdots, m; l = 1, \cdots, d$ in Banach spaces and the sequences in $(2.4)$ converges to the fixed points of $ U_{k,l,n}, P_{k,n}, n \geq 0, k = 1, \cdots, m; l = 1, \cdots, d$ respectively, which are the solutions, $ u_{k,l} = \lim \limits_{n \to \infty} u_{k,l,n}, p_{k} = \lim \limits_{n \to \infty} p_{k,n}, k = 1, \cdots, m; l = 1, \cdots, d$, of the system $(2.1)$.

\section{Variational iteration method for the system of incompressible fluid flow} 
In this section, we apply our technique of generalized He's variational iteration method to create an iterative scheme that can guarantee existence of unique solution for the system of incompressible fluid flow which is modeled by Navier-Stokes Equations, see [21,22,3,20,17]: given $\textit{\textbf{f}},\textit{\textbf{g}}$ and time $T > 0$, to find $\textit{\textbf{u}} = (u_{1},u_{2},u_{3}): \Omega \times \left[ 0,T\right] \to \mathbb{R}^d$ in $L^{\infty} \left( \left( 0,T \right); \left( L^{2} \left( \Omega \right) \right)^{n} \right) \bigcap L^{2} \left( \left( 0,T \right); \left(  H^{1}_{0} \left( \Omega \right) \right)^{n} \right) $ and $p : \Omega   \times \left[ 0,T\right] \to \mathbb{R}$ in $ L^{\frac{n+2}{n+1}} \left( \left( 0,T \right); L^{\frac{n+2}{n+1}} \left( \Omega \right) \right) $such that
\begin{align} \left\{ \begin{array}{ll}
\dfrac{\partial  \textit{\textbf{u}}}{\partial t} + \left(  \textit{\textbf{u}}.\nabla \right)  \textit{\textbf{u}} - \mathrm{v} \triangle  \textit{\textbf{u}} = -\nabla p + \textit{\textbf{f}} \ &  \mbox{ $ \mathrm{in} \ \Omega \times \left( 0,T\right]$} , \\
\nabla. \textit{\textbf{u}} = 0  & \mbox{ $\mathrm{in} \ \Omega \times \left( 0,T\right]$}, \\
\textit{\textbf{u}} =  \textit{\textbf{g}}  &  \mbox{ $\mathrm{at} \ \partial \Omega \times \left( 0,T\right]$},  \\
\textit{\textbf{u}}(\textit{\textbf{x}},0) =  \textit{\textbf{u}}_{0}(\textit{\textbf{x}}) \  &  \mbox{ $\mathrm{in} \ \Omega$}. 
\end{array} 
\right.
\end{align} 
where $ \textit{\textbf{u}}$ is the velocity field, $p$ is the pressure field, $\mathrm{v} (> 0)$ is the coefficient of kinematical viscosity and $\textit{\textbf{f}}$ is the body force, we normalize $(3.1_{1})$ by divide both sides on $ - \mathrm{v}$ to reform the system to 
\begin{align} \left\{ \begin{array}{ll}
\triangle \textit{\textbf{u}}  - \frac{1}{\mathrm{v}} \left( \dfrac{\partial \textit{\textbf{u}}}{\partial t} + \left(  \textit{\textbf{u}}.\nabla \right)  \textit{\textbf{u}} \right)   = \frac{1}{\mathrm{v}} \left( \nabla p - \textit{\textbf{f}} \right) \ &  \mbox{ $ \mathrm{in} \ \Omega \times \left( 0,T\right]$} , \\
\nabla. \textit{\textbf{u}} = 0  & \mbox{ $\mathrm{in} \ \Omega \times \left( 0,T\right]$}, \\
 \textit{\textbf{u}} =  \textit{\textbf{g}}  &  \mbox{ $\mathrm{at} \ \partial \Omega \times \left( 0,T\right]$},  \\
 \textit{\textbf{u}}(\textit{\textbf{x}},0) =  \textit{\textbf{u}}_{0}(\textit{\textbf{x}}) \  &  \mbox{ $\mathrm{in} \ \Omega$}. \\
\end{array} 
\right.
\end{align} 
Since the unknowns $\textit{\textbf{u}},p$ do not appear in $(3.1)$ in a symmetric way, because the pressure plays the role of reaction force associated with the isochoricity constraint $\nabla .  \textit{\textbf{u}} = 0$, see [5], we take $\nabla$ on both sides of $(3.1_{1})$ to get the field $p$ as a solution of the following Neumann problem
\begin{align} \left\{ \begin{array}{ll}
\triangle p =  - \nabla \left( \left(  \textit{\textbf{u}}.\nabla \right)  \textit{\textbf{u}} - \textit{\textbf{f}}  \right) &  \mbox{ $ \mathrm{in} \ \Omega \times \left( 0,T\right]$} , \\
\dfrac{\partial p}{\partial n} = - \left( \mathrm{v} \triangle  \textit{\textbf{u}} + \textit{\textbf{f}} \right). \textit{\textbf{n}}  & \mbox{ $\mathrm{at} \ \partial \Omega \times \left( 0,T\right]$}, \\
\end{array} 
\right.
\end{align} 
where $\textit{\textbf{n}}$ is the unit outer normal to $\partial \Omega$. $(3.2_{1})$ and $(3.3_{1})$ can be expanded in components to the following:
\begin{align*} 
\frac{\partial^{2} u_{1}}{\partial x_{1}^{2}} + \frac{\partial^{2} u_{1}}{\partial x_{2}^{2}} + \frac{\partial^{2} u_{1}} {\partial x_{3}^{2}}  &- \frac{1}{\mathrm{v}} \left( \frac{\partial u_{1}}{\partial t} + u_{1}\frac{\partial u_{1}}{\partial x_{1}} + u_{2}\frac{\partial u_{1}}{\partial x_{2}} + u_{3}\frac{\partial u_{1}}{\partial x_{3}} \right)  =  \frac{1}{\mathrm{v}} \left( \frac{\partial p}{\partial x_{1}} - f_{1} \right) 
\\
\frac{\partial^{2} u_{2}}{\partial x_{1}^{2}} + \frac{\partial^{2} u_{2}}{\partial x_{2}^{2}} + \frac{\partial^{2} u_{2}}{\partial x_{3}^{2}} &- \frac{1}{\mathrm{v}} \left( \frac{\partial u_{2}}{\partial t} + u_{1}\frac{\partial u_{2}}{\partial x_{1}} + u_{2}\frac{\partial u_{2}}{\partial x_{2}} + u_{3}\frac{\partial u_{2}}{\partial x_{3}} \right) =  \frac{1}{\mathrm{v}} \left( \frac{\partial p}{\partial x_{2}} - f_{2} \right) 
 \\
\frac{\partial^{2} u_{3}}{\partial x_{1}^{2}} + \frac{\partial^{2} u_{3}}{\partial x_{2}^{2}} + \frac{\partial^{2} u_{3}}{\partial x_{3}^{2}} &- \frac{1}{\mathrm{v}} \left( \frac{\partial u_{3}}{\partial t} + u_{1}\frac{\partial u_{3}}{\partial x_{1}} + u_{2}\frac{\partial u_{3}}{\partial x_{2}} + u_{3}\frac{\partial u_{3}}{\partial x_{3}} \right)  =  \frac{1}{\mathrm{v}} \left( \frac{\partial p}{\partial x_{3}} - f_{3} \right) 
\\
\frac{\partial^{2} p}{\partial x_{1}^{2}} + \frac{\partial^{2} p}{\partial x_{2}^{2}} + \frac{\partial^{2} p} {\partial x_{3}^{2}} &= \frac{\partial}{\partial x_{1}} \left( u_{1} \frac{\partial u_{1}}{\partial x_{1}} + u_{2} \frac{\partial u_{1}}{\partial x_{2}} + u_{3} \frac{\partial u_{1}}{\partial x_{3}} - f_{1} \right) + \frac{\partial}{\partial x_{2}} \left( u_{1} \frac{\partial u_{2}}{\partial x_{1}} + u_{2} \frac{\partial u_{2}}{\partial x_{2}} + u_{3} \frac{\partial u_{2}}{\partial x_{3}} - f_{2} \right) \\
&+ \frac{\partial}{\partial x_{3}} \left( u_{1} \frac{\partial u_{3}}{\partial x_{1}} + u_{2} \frac{\partial u_{3}}{\partial x_{2}} + u_{3} \frac{\partial u_{3}}{\partial x_{3}} - f_{3} \right). 
\end{align*}
Then, the correctional functionals for above equations in $x_{1}-$ direction are    
\begin{align*}
u_{1,n+1}(x_{1},x_{2},x_{3},t) &= u_{1,n}(x_{1},x_{2},x_{3},t) + \int_{\Omega_{x_{1}}} \lambda_{1}(X_{1}) \bigg[ \frac{\partial^{2} u_{1,n}(X_{1},x_{2},x_{3},t)}{\partial X_{1}^{2}} + \frac{\partial^{2} \tilde{u}_{1,n}(X_{1},x_{2},x_{3},t)}{\partial x_{2}^{2}} \\
&+ \frac{\partial^{2} \tilde{u}_{1,n}(X_{1},x_{2},x_{3},t)} {\partial x_{3}^{2}} - \frac{1}{\mathrm{v}} \bigg( \frac{\partial \tilde{u}_{1,n}(X_{1},x_{2},x_{3},t)}{\partial t} + \tilde{u}_{1,n}(X_{1},x_{2},x_{3},t) \frac{\partial \tilde{u}_{1,n}(X_{1},x_{2},x_{3},t)}{\partial X_{1}} \\
&+ \tilde{u}_{2,n}(X_{1},x_{2},x_{3},t) \frac{\partial \tilde{u}_{1,n}(X_{1},x_{2},x_{3},t)}{\partial x_{2}} + \tilde{u}_{3,n}(X_{1},x_{2},x_{3},t) \frac{\partial \tilde{u}_{1,n}(X_{1},x_{2},x_{3},t)}{\partial x_{3}} \bigg) \\
&- \frac{1}{\mathrm{v}} \left( \frac{\partial \tilde{p}_{n}(X_{1},x_{2},x_{3},t)}{\partial X_{1}} - f_{1} (X_{1},x_{2},x_{3},t) \right) 
\bigg] d X_{1},
\\
u_{2,n+1}(x_{1},x_{2},x_{3},t) &= u_{2,n}(x_{1},x_{2},x_{3},t) + \int_{\Omega_{x_{1}}} \lambda_{2}(X_{1}) \bigg[ \frac{\partial^{2} u_{2,n}(X_{1},x_{2},x_{3},t)}{\partial X_{1}^{2}} + \frac{\partial^{2} \tilde{u}_{2,n}(X_{1},x_{2},x_{3},t)}{\partial x_{2}^{2}} \\
&+ \frac{\partial^{2} \tilde{u}_{2,n}(X_{1},x_{2},x_{3},t)} {\partial x_{3}^{2}} - \frac{1}{\mathrm{v}} \bigg( \frac{\partial \tilde{u}_{2,n}(X_{1},x_{2},x_{3},t)}{\partial t} + \tilde{u}_{1,n}(X_{1},x_{2},x_{3},t) \frac{\partial \tilde{u}_{2,n}(X_{1},x_{2},x_{3},t)}{\partial X_{1}} \\
&+ \tilde{u}_{2,n}(X_{1},x_{2},x_{3},t) \frac{\partial \tilde{u}_{2,n}(X_{1},x_{2},x_{3},t)}{\partial x_{2}} + \tilde{u}_{3,n}(X_{1},x_{2},x_{3},t) \frac{\partial \tilde{u}_{2,n}(X_{1},x_{2},x_{3},t)}{\partial x_{3}} \bigg) \\
&- \frac{1}{\mathrm{v}} \left( \frac{\partial \tilde{p}_{n}(X_{1},x_{2},x_{3},t)}{\partial x_{2}} - f_{2} (X_{1},x_{2},x_{3},t) \right) 
\bigg] d X_{1},
\\ 
u_{3,n+1}(x_{1},x_{2},x_{3},t) &= u_{3,n}(x_{1},x_{2},x_{3},t) + \int_{\Omega_{x_{1}}} \lambda_{3}(X_{1}) \bigg[ \frac{\partial^{2} u_{3,n}(X_{1},x_{2},x_{3},t)}{\partial X_{1}^{2}} + \frac{\partial^{2} \tilde{u}_{3,n}(X_{1},x_{2},x_{3},t)}{\partial x_{2}^{2}} \\
&+ \frac{\partial^{2} \tilde{u}_{3,n}(X_{1},x_{2},x_{3},t)} {\partial x_{3}^{2}} - \frac{1}{\mathrm{v}} \bigg( \frac{\partial \tilde{u}_{3,n}(X_{1},x_{2},x_{3},t)}{\partial t} + \tilde{u}_{1,n}(X_{1},x_{2},x_{3},t) \frac{\partial \tilde{u}_{3,n}(X_{1},x_{2},x_{3},t)}{\partial X_{1}} \\
&+ \tilde{u}_{2,n}(X_{1},x_{2},x_{3},t) \frac{\partial \tilde{u}_{3,n}(X_{1},x_{2},x_{3},t)}{\partial x_{2}} + \tilde{u}_{3,n}(X_{1},x_{2},x_{3},t) \frac{\partial \tilde{u}_{3,n}(X_{1},x_{2},x_{3},t)}{\partial x_{3}} \bigg) \\
&- \frac{1}{\mathrm{v}} \left( \frac{\partial \tilde{p}_{n}(X_{1},x_{2},x_{3},t)}{\partial x_{3}} - f_{3} (X_{1},x_{2},x_{3},t) \right) 
\bigg] d X_{1}, 
\\
p_{n+1}(x_{1},x_{2},x_{3},t) &= p_{n}(x_{1},x_{2},x_{3},t) + \int_{\Omega_{x_{1}}} \mu(X_{1}) \bigg[ \frac{\partial^{2} p_{n}(X_{1},x_{2},x_{3},t)}{\partial X_{1}^{2}} + \frac{\partial^{2} \tilde{p}_{n}(X_{1},x_{2},x_{3},t)}{\partial x_{2}^{2}} + \frac{\partial^{2} \tilde{p}_{n}(X_{1},x_{2},x_{3},t)} {\partial x_{3}^{2}} \\
&+ \frac{\partial}{\partial X_{1}} \bigg( \tilde{u}_{1,n}(X_{1},x_{2},x_{3},t) \frac{\partial \tilde{u}_{1,n}(X_{1},x_{2},x_{3},t)}{\partial X_{1}} + \tilde{u}_{2,n}(X_{1},x_{2},x_{3},t) \frac{\partial \tilde{u}_{1,n}(X_{1},x_{2},x_{3},t)}{\partial x_{2}} \\
&+ \tilde{u}_{3,n}(X_{1},x_{2},x_{3},t) \frac{\partial \tilde{u}_{1,n}(X_{1},x_{2},x_{3},t)}{\partial x_{3}} - f_{1}(X_{1},x_{2},x_{3},t) \bigg) \\
&+ \frac{\partial}{\partial x_{2}} \bigg( \tilde{u}_{1,n}(X_{1},x_{2},x_{3},t) \frac{\partial \tilde{u}_{2,n}(X_{1},x_{2},x_{3},t)}{\partial X_{1}} + \tilde{u}_{2,n}(X_{1},x_{2},x_{3},t) \frac{\partial \tilde{u}_{2,n}(X_{1},x_{2},x_{3},t)}{\partial x_{2}} \\
&+ \tilde{u}_{3,n}(X_{1},x_{2},x_{3},t) \frac{\partial \tilde{u}_{2,n}(X_{1},x_{2},x_{3},t)}{\partial x_{3}} - f_{2}(X_{1},x_{2},x_{3},t) \bigg) \\
&+ \frac{\partial}{\partial x_{3}} \bigg( \tilde{u}_{1,n}(X_{1},x_{2},x_{3},t) \frac{\partial \tilde{u}_{3,n}(X_{1},x_{2},x_{3},t)}{\partial X_{1}} + \tilde{u}_{2,n}(X_{1},x_{2},x_{3},t) \frac{\partial \tilde{u}_{3,n}(X_{1},x_{2},x_{3},t)}{\partial x_{2}} \\
&+ \tilde{u}_{3,n}(X_{1},x_{2},x_{3},t) \frac{\partial \tilde{u}_{3,n}(X_{1},x_{2},x_{3},t)}{\partial x_{3}} - f_{3}(X_{1},x_{2},x_{3},t) \bigg) \bigg] d X_{1}.
\end{align*}
Making the above correctional functionals stationary,  
\begin{align*}
\delta u_{1,n+1}(x_{1},x_{2},x_{3},t) &= \delta u_{1,n}(x_{1},x_{2},x_{3},t) + \left. \frac{\partial \delta u_{1,n}(X_{1},x_{2},x_{3},t)}{\partial X_{1}} \lambda_{1}(X_{1}) \right|_{X_{1} = x_{1}} - \left. \delta u_{1,n}(X_{1},x_{2},x_{3},t) \frac{d \lambda_{1}(X_{1})}{d X_{1}}  \right|_{X_{1} = x_{1}} \\
&+ \int_{\Omega_{x_{1}}} \delta u_{1,n}(X_{1},x_{2},x_{3},t) \frac{d^{2} \lambda_{1}(X_{1})}{d X^{2}_{1}} dX_{1} = 0,
\\
\delta u_{2,n+1}(x_{1},x_{2},x_{3},t) &= \delta u_{2,n}(x_{1},x_{2},x_{3},t) + \left. \frac{\partial \delta u_{2,n}(X_{1},x_{2},x_{3},t)}{\partial X_{1}} \lambda_{2}(X_{1}) \right|_{X_{1} = x_{1}} - \left. \delta u_{2,n}(X_{1},x_{2},x_{3},t) \frac{d \lambda_{2}(X_{1})}{d X_{1}}  \right|_{X_{1} = x_{1}} \\
&+ \int_{\Omega_{x_{1}}} \delta u_{2,n}(X_{1},x_{2},x_{3},t) \frac{d^{2} \lambda_{2}(X_{1})}{d X^{2}_{1}} dX_{1} = 0,
\\
\delta u_{3,n+1}(x_{1},x_{2},x_{3},t) &= \delta u_{3,n}(x_{1},x_{2},x_{3},t) + \left. \frac{\partial \delta u_{3,n}(X_{1},x_{2},x_{3},t)}{\partial X_{1}} \lambda_{3}(X_{1}) \right|_{X_{1} = x_{1}} - \left. \delta u_{3,n}(X_{1},x_{2},x_{3},t) \frac{d \lambda_{3}(X_{1})}{d X_{1}}  \right|_{X_{1} = x_{1}} \\
&+ \int_{\Omega_{x_{1}}} \delta u_{3,n}(X_{1},x_{2},x_{3},t) \frac{d^{2} \lambda_{3}(X_{1})}{d X^{2}_{1}} dX_{1} = 0, 
\\
\delta p_{n+1}(x_{1},x_{2},x_{3},t) &= \delta p_{n}(x_{1},x_{2},x_{3},t) + \left. \frac{\partial \delta p_{n}(X_{1},x_{2},x_{3},t)}{\partial X_{1}} \mu(X_{1}) \right|_{X_{1} = x_{1}} - \left. \delta p_{n}(X_{1},x_{2},x_{3},t) \frac{d \mu(X_{1})}{d X_{1}}  \right|_{X_{1} = x_{1}} \\
&+ \int_{\Omega_{x_{1}}} \delta p_{n}(X_{1},x_{2},x_{3},t) \frac{d^{2} \mu(X_{1})}{d X^{2}_{1}} dX_{1} = 0
\end{align*}
yields the following stationary conditions: for $i = 1,2,3$, that
\begin{align*}
\delta u_{i,n}: \frac{d^{2} \lambda_{i}}{d X^{2}_{1}} &= 0  & \ \ \ \delta p_{n}: \frac{d^{2} \mu}{d X^{2}_{1}} = 0 \\
\frac{\partial \delta u_{i,n}}{\partial X_{1}}: \left. \mu(X_{1}) \right|_{X_{1} = x_{1}} &= 1  \ \ \ \ \ \ \ \ \ \ \ \ \ \ \ \ \ \mathrm{and}& \ \ \ \frac{\partial \delta p_{n}}{\partial X_{1}}: \left. \mu(X_{1}) \right|_{X_{1} = x_{1}} = 1 \\
\delta u_{i,n}: \left. \frac{d \lambda_{i} (X_{1})}{d X_{1}} \right|_{X_{1} = x_{1}} &= 0  & \ \ \ \delta p_{n}: \left. \frac{d \mu (X_{1})}{d X_{1}} \right|_{X_{1} = x_{1}} = 0.
\end{align*}
Then, for $i = 1,2,3$, the Lagrange multipliers are $\lambda_{i} = 1$ and $\mu = 1$; and the desired iterative scheme is
\begin{align*}
u_{1,n+1}(x_{1},x_{2},x_{3},t) &= u_{1,n}(x_{1},x_{2},x_{3},t) + \int_{\Omega_{x_{1}}} \bigg[ \frac{\partial^{2} u_{1,n}(X_{1},x_{2},x_{3},t)}{\partial X_{1}^{2}} + \frac{\partial^{2} u_{1,n}(X_{1},x_{2},x_{3},t)}{\partial x_{2}^{2}} \\
&+ \frac{\partial^{2} u_{1,n}(X_{1},x_{2},x_{3},t)} {\partial x_{3}^{2}} - \frac{1}{\mathrm{v}} \bigg( \frac{\partial u_{1,n}(X_{1},x_{2},x_{3},t)}{\partial t} + u_{1,n}(X_{1},x_{2},x_{3},t) \frac{\partial u_{1,n}(X_{1},x_{2},x_{3},t)}{\partial X_{1}} \\
&+ u_{2,n}(X_{1},x_{2},x_{3},t) \frac{\partial u_{1,n}(X_{1},x_{2},x_{3},t)}{\partial x_{2}} + u_{3,n}(X_{1},x_{2},x_{3},t) \frac{\partial u_{1,n}(X_{1},x_{2},x_{3},t)}{\partial x_{3}} \bigg) \\
&- \frac{1}{\mathrm{v}} \left( \frac{\partial p_{n}(X_{1},x_{2},x_{3},t)}{\partial X_{1}} - f_{1} (X_{1},x_{2},x_{3},t) \right) \bigg] d X_{1},
\end{align*} 
\begin{align*}
u_{2,n+1}(x_{1},x_{2},x_{3},t) &= u_{2,n}(x_{1},x_{2},x_{3},t) + \int_{\Omega_{x_{1}}} \bigg[ \frac{\partial^{2} u_{2,n}(X_{1},x_{2},x_{3},t)}{\partial X_{1}^{2}} + \frac{\partial^{2} u_{2,n}(X_{1},x_{2},x_{3},t)}{\partial x_{2}^{2}} \\
&+ \frac{\partial^{2} u_{2,n}(X_{1},x_{2},x_{3},t)} {\partial x_{3}^{2}} - \frac{1}{\mathrm{v}} \bigg( \frac{\partial u_{2,n}(X_{1},x_{2},x_{3},t)}{\partial t} + u_{1,n}(X_{1},x_{2},x_{3},t) \frac{\partial u_{2,n}(X_{1},x_{2},x_{3},t)}{\partial X_{1}} \\
&+ u_{2,n}(X_{1},x_{2},x_{3},t) \frac{\partial u_{2,n}(X_{1},x_{2},x_{3},t)}{\partial x_{2}} + u_{3,n}(X_{1},x_{2},x_{3},t) \frac{\partial u_{2,n}(X_{1},x_{2},x_{3},t)}{\partial x_{3}} \bigg) \\
&- \frac{1}{\mathrm{v}} \left( \frac{\partial p_{n}(X_{1},x_{2},x_{3},t)}{\partial x_{2}} - f_{2} (X_{1},x_{2},x_{3},t) \right) \bigg] d X_{1},
\end{align*} 
\begin{align*} 
u_{3,n+1}(x_{1},x_{2},x_{3},t) &= u_{3,n}(x_{1},x_{2},x_{3},t) + \int_{\Omega_{x_{1}}} \bigg[ \frac{\partial^{2} u_{3,n}(X_{1},x_{2},x_{3},t)}{\partial X_{1}^{2}} + \frac{\partial^{2} u_{3,n}(X_{1},x_{2},x_{3},t)}{\partial x_{2}^{2}} \\
&+ \frac{\partial^{2} u_{3,n}(X_{1},x_{2},x_{3},t)} {\partial x_{3}^{2}} - \frac{1}{\mathrm{v}} \bigg( \frac{\partial u_{3,n}(X_{1},x_{2},x_{3},t)}{\partial t} + u_{1,n}(X_{1},x_{2},x_{3},t) \frac{\partial u_{3,n}(X_{1},x_{2},x_{3},t)}{\partial X_{1}} \\
&+ u_{2,n}(X_{1},x_{2},x_{3},t) \frac{\partial u_{3,n}(X_{1},x_{2},x_{3},t)}{\partial x_{2}} + u_{3,n}(X_{1},x_{2},x_{3},t) \frac{\partial u_{3,n}(X_{1},x_{2},x_{3},t)}{\partial x_{3}} \bigg) \\
&- \frac{1}{\mathrm{v}} \left( \frac{\partial p_{n}(X_{1},x_{2},x_{3},t)}{\partial x_{3}} - f_{3} (X_{1},x_{2},x_{3},t) \right) \bigg] d X_{1}, 
\end{align*} 
\begin{align*}
p_{n+1}(x_{1},x_{2},x_{3},t) &= p_{n}(x_{1},x_{2},x_{3},t) + \int_{\Omega_{x_{1}}} \bigg[ \frac{\partial^{2} p_{n}(X_{1},x_{2},x_{3},t)}{\partial X_{1}^{2}} + \frac{\partial^{2} p_{n}(X_{1},x_{2},x_{3},t)}{\partial x_{2}^{2}} + \frac{\partial^{2} p_{n}(X_{1},x_{2},x_{3},t)} {\partial x_{3}^{2}} \\
&+ \frac{\partial}{\partial X_{1}} \bigg( u_{1,n}(X_{1},x_{2},x_{3},t) \frac{\partial u_{1,n}(X_{1},x_{2},x_{3},t)}{\partial X_{1}} + u_{2,n}(X_{1},x_{2},x_{3},t) \frac{\partial u_{1,n}(X_{1},x_{2},x_{3},t)}{\partial x_{2}} \\
&+ u_{3,n}(X_{1},x_{2},x_{3},t) \frac{\partial u_{1,n}(X_{1},x_{2},x_{3},t)}{\partial x_{3}} - f_{1}(X_{1},x_{2},x_{3},t) \bigg) \\
&+ \frac{\partial}{\partial x_{2}} \bigg( u_{1,n}(X_{1},x_{2},x_{3},t) \frac{\partial u_{2,n}(X_{1},x_{2},x_{3},t)}{\partial X_{1}} + u_{2,n}(X_{1},x_{2},x_{3},t) \frac{\partial u_{2,n}(X_{1},x_{2},x_{3},t)}{\partial x_{2}} \\
&+ u_{3,n}(X_{1},x_{2},x_{3},t) \frac{\partial u_{2,n}(X_{1},x_{2},x_{3},t)}{\partial x_{3}} - f_{2}(X_{1},x_{2},x_{3},t) \bigg) \\
&+ \frac{\partial}{\partial x_{3}} \bigg( u_{1,n}(X_{1},x_{2},x_{3},t) \frac{\partial u_{3,n}(X_{1},x_{2},x_{3},t)}{\partial X_{1}} + u_{2,n}(X_{1},x_{2},x_{3},t) \frac{\partial u_{3,n}(X_{1},x_{2},x_{3},t)}{\partial x_{2}} \\
&+ u_{3,n}(X_{1},x_{2},x_{3},t) \frac{\partial u_{3,n}(X_{1},x_{2},x_{3},t)}{\partial x_{3}} - f_{3}(X_{1},x_{2},x_{3},t) \bigg) \bigg] d X_{1}.
\end{align*} 

\section{Variational iteration method for the system of incompressible micropolar fluid flow} 
Another application of the generalized He's variational iteration method is to create an iterative scheme that can guarantee existence of unique solution for the system of incompressible micropolar fluid flow which is modeled by Navier-Stokes equations and microrotational velocity equations, see[4,18,19]: given $\textit{\textbf{f}}_{1},\textit{\textbf{f}}_{2},\textit{\textbf{g}},\textit{\textbf{q}}$ and time $T > 0$, to find $\textit{\textbf{u}} = (u_{1},u_{2},u_{3}): \Omega \times \left[ 0,T\right] \to \mathbb{R}^{3}$ in $L^{\infty} \left( \left( 0,T \right); L^{2} \left( \Omega \right) \right) \bigcap L^{2} \left( \left( 0,T \right); V_{g} \right)$, $\textit{\textbf{w}} = (w_{1},w_{2},w_{3}): \Omega \times \left[ 0,T\right] \to \mathbb{R}^{3}$ in $L^{\infty} \left( \left( 0,T \right); L^{2} \left( \Omega \right) \right) \bigcap L^{2} \left( \left( 0,T \right); H^{1}_{q} \left( \Omega \right) \right) $ and $p : \Omega   \times \left[ 0,T\right] \to \mathbb{R}$ in $ L^{2} \left( \left( 0,T \right); L^{2}_{0} \left( \Omega \right) \right)$ such that
\begin{align} \left\{ \begin{array}{lc} 
\dfrac{\partial \textit{\textbf{u}}}{\partial t} - \left( \mathrm{v} + \mathrm{v_{r}} \right)  \triangle \textit{\textbf{u}} + \left( \textit{\textbf{u}}.\nabla \right) \textit{\textbf{u}} + \nabla p   = 2 \mathrm{v_{r}} \nabla \times \textit{\textbf{w}} + \textit{\textbf{f}}_{1}  &  \mbox{ $ \mathrm{in} \ \Omega \times \left( 0,T\right]$} , \\
\dfrac{\partial \textit{\textbf{w}}}{\partial t} - \left( \mathrm{c_{a}} + \mathrm{c_{d}} \right)  \triangle \textit{\textbf{w}} - \left( \mathrm{c_{0}} + \mathrm{c_{d}} - \mathrm{c_{a}} \right) \nabla \left( \nabla.\textit{\textbf{w}} \right) + \left( \textit{\textbf{u}}.\nabla \right) \textit{\textbf{w}} + 4\mathrm{v_{r}} \textit{\textbf{w}}
 = 2 \mathrm{v_{r}} \nabla \times \textit{\textbf{u}} + \textit{\textbf{f}}_{2}  &  \mbox{ $ \mathrm{in} \ \Omega \times \left( 0,T\right]$} , \\
\nabla.\textit{\textbf{u}} = 0  & \mbox{ $\mathrm{in} \ \Omega \times \left( 0,T\right]$}, \\
\textit{\textbf{u}} = \textit{\textbf{g}}, \textit{\textbf{w}} = \textit{\textbf{q}}   &  \mbox{ $\mathrm{on} \ \partial \Omega \times \left( 0,T\right]$},  \\
\textit{\textbf{u}}(\textit{\textbf{x}},0) = \textit{\textbf{u}}_{0}(\textit{\textbf{x}}), \textit{\textbf{w}}\nabla (\textit{\textbf{x}},0) = \textit{\textbf{w}}_{0}(\textit{\textbf{x}})  &  \mbox{ $\mathrm{in} \ \Omega$}. \\
\end{array} 
\right.
\end{align} 
where $ V_{g} = \left\lbrace  \textit{\textbf{v}} \in H^{1} (\Omega): \left. \textit{\textbf{v}} \right| _{\partial \Omega} = \textit{\textbf{g}}, \nabla . \textit{\textbf{v}} = 0 \ \mathrm{in} \ \Omega \right\rbrace, H^{1}_{q} \left( \Omega \right) = \left\lbrace  \textit{\textbf{v}} \in H^{1} (\Omega): \left. \textit{\textbf{v}} \right| _{\partial \Omega} = \textit{\textbf{q}} \right\rbrace $ and $L^{2}_{0} \left( \Omega \right) = \{ p \in L^{2} \left( \Omega \right): \ \int_{\Omega} p d\Omega = 0 \}$. $\textit{\textbf{u}}$ is the fluid velocity, $\textit{\textbf{w}}$ the microrotation field(the angular velocity field of rotation of particles) and $p$ the fluid kinematic pressure. The fields $\textit{\textbf{f}}_{1}$ and $\textit{\textbf{f}}_{2}$ are the external body force and moment(torgue) respectively. The positive constants $\mathrm{v},\mathrm{v_{r}}, \mathrm{c_{a}}$ and $\mathrm{c_{d}}$ represent viscosity coefficients, $\mathrm{v}$ is the Newtonian viscosity and $\mathrm{v_{r}}$ is the microrotation viscosity. The constants $\mathrm{c_{0}},\mathrm{c_{d}}$ and $\mathrm{c_{a}}$ satisfy $\mathrm{c_{0}}+\mathrm{c_{d}} > \mathrm{c_{a}}$. 
 We normalize $(4.1_{1})$ and $(4.1_{2})$ by divide both sides by $- \left( \mathrm{v} + \mathrm{v_{r}} \right)$ and $- \left( \mathrm{c_{a}} + \mathrm{c_{d}} \right)$ respectively, to reform the system $(4.1)$ to
\begin{align} \left\{ \begin{array}{lc} 
\triangle \textit{\textbf{u}} - \frac{1}{\left( \mathrm{v} + \mathrm{v_{r}} \right)} \left( \dfrac{\partial \textit{\textbf{u}}}{\partial t} + \left( \textit{\textbf{u}}.\nabla \right) \textit{\textbf{u}} + \nabla p \right)   = - \frac{1}{\left( \mathrm{v} + \mathrm{v_{r}} \right)} \left( 2 \mathrm{v_{r}} \nabla \times \textit{\textbf{w}} + \textit{\textbf{f}}_{1} \right) &  \mbox{ $ \mathrm{in} \ \Omega \times \left( 0,T\right]$} ,
\\
\triangle \textit{\textbf{w}} - \frac{1}{\left( \mathrm{c_{a}} + \mathrm{c_{d}} \right)} \left( \dfrac{\partial \textit{\textbf{w}}}{\partial t}  - \left( \mathrm{c_{0}} + \mathrm{c_{d}} - \mathrm{c_{a}} \right) \nabla \left( \nabla.\textit{\textbf{w}} \right) + \left( \textit{\textbf{u}}.\nabla \right) \textit{\textbf{w}} + 4\mathrm{v_{r}} \textit{\textbf{w}} \right)
 = - \frac{1}{\left( \mathrm{c_{a}} + \mathrm{c_{d}} \right)} \left( 2 \mathrm{v_{r}} \nabla \times \textit{\textbf{u}} + \textit{\textbf{f}}_{2} \right)  &  \mbox{ $ \mathrm{in} \ \Omega \times \left( 0,T\right]$} , \\
\nabla.\textit{\textbf{u}} = 0  & \mbox{ $\mathrm{in} \ \Omega \times \left( 0,T\right]$}, \\
\textit{\textbf{u}} = \textit{\textbf{g}}, \textit{\textbf{w}} = \textit{\textbf{q}}   &  \mbox{ $\mathrm{on} \ \partial \Omega \times \left( 0,T\right]$},  \\
\textit{\textbf{u}}(\textit{\textbf{x}},0) = \textit{\textbf{u}}_{0}(\textit{\textbf{x}}), \textit{\textbf{w}}(\textit{\textbf{x}},0) = \textit{\textbf{w}}_{0}(\textit{\textbf{x}})  &  \mbox{ $\mathrm{in} \ \Omega$}. \\
\end{array} 
\right. 
\end{align} 
Since the unknowns $\textit{\textbf{u}},\textit{\textbf{w}},p$ do not appear in symmetric way in $(4.1_{1})$, the field $p$ can be obtained by taking $\nabla$ on both sides of $(4.1_{1})$ as a solution of the following Neumann problem
\begin{align} \left\{ \begin{array}{ll}
\triangle p =  - \nabla \left( \left( \textit{\textbf{u}}.\nabla \right) \textit{\textbf{u}} - 2 \mathrm{v_{r}} \nabla \times \textit{\textbf{w}} - \textit{\textbf{f}}_{1}  \right) &  \mbox{ $ \mathrm{in} \ \Omega \times \left( 0,T\right]$} , \\
\dfrac{\partial p}{\partial n} = - \left( \mathrm{v} \triangle \textit{\textbf{u}} + 2 \mathrm{v_{r}} \nabla \times \textit{\textbf{u}} + \textit{\textbf{f}}_{1}  \right). \textit{\textbf{n}}  & \mbox{ $\mathrm{at} \ \partial \Omega \times \left( 0,T\right]$}, \\
\end{array} 
\right.
\end{align}  
where $\textit{\textbf{n}}$ is the unit outer normal to $\partial \Omega$. $(4.2_{1})$ and $(4.3_{1})$ can be expanded in components to the following:
\begin{align*} 
\frac{\partial^{2} u_{1}}{\partial x_{1}^{2}} + \frac{\partial^{2} u_{1}}{\partial x_{2}^{2}} + \frac{\partial^{2} u_{1}} {\partial x_{3}^{2}} - \frac{1}{\left( \mathrm{v} + \mathrm{v_{r}} \right)} \left( \frac{\partial u_{1}}{\partial t} + u_{1}\frac{\partial u_{1}}{\partial x_{1}} + u_{2}\frac{\partial u_{1}}{\partial x_{2}} + u_{3}\frac{\partial u_{1}}{\partial x_{3}} + \frac{\partial p}{\partial x_{1}}  \right)  &= - \frac{1}{\left( \mathrm{v} + \mathrm{v_{r}} \right)} \left( 2 \mathrm{v_{r}} \left( \frac{\partial w_{3}}{\partial x_{2}} - \frac{\partial w_{2}}{\partial x_{3}} \right) + f_{1,1} \right), 
\end{align*}  
\begin{align*} 
\frac{\partial^{2} u_{2}}{\partial x_{1}^{2}} + \frac{\partial^{2} u_{2}}{\partial x_{2}^{2}} + \frac{\partial^{2} u_{2}} {\partial x_{3}^{2}}  - \frac{1}{\left( \mathrm{v} + \mathrm{v_{r}} \right)} \left( \frac{\partial u_{2}}{\partial t} + u_{1}\frac{\partial u_{2}}{\partial x_{1}} + u_{2}\frac{\partial u_{2}}{\partial x_{2}} + u_{3}\frac{\partial u_{2}}{\partial x_{3}} + \frac{\partial p}{\partial x_{2}} \right)  &= - \frac{1}{\left( \mathrm{v} + \mathrm{v_{r}} \right)} \left( 2 \mathrm{v_{r}} \left( \frac{\partial w_{1}}{\partial x_{3}} - \frac{\partial w_{3}}{\partial x_{1}} \right) + f_{1,2} \right), 
\end{align*}  
\begin{align*} 
\frac{\partial^{2} u_{1}}{\partial x_{1}^{2}} + \frac{\partial^{2} u_{1}}{\partial x_{2}^{2}} + \frac{\partial^{2} u_{3}} {\partial x_{3}^{2}} - \frac{1}{\left( \mathrm{v} + \mathrm{v_{r}} \right)} \left( \frac{\partial u_{3}}{\partial t} + u_{1}\frac{\partial u_{3}}{\partial x_{1}} + u_{2}\frac{\partial u_{3}}{\partial x_{2}} + u_{3} \frac{\partial u_{3}}{\partial x_{3}} + \frac{\partial p}{\partial x_{3}} \right)  &= - \frac{1}{\left( \mathrm{v} + \mathrm{v_{r}} \right)} \left( 2 \mathrm{v_{r}} \left( \frac{\partial w_{2}}{\partial x_{1}} - \frac{\partial w_{1}}{\partial x_{2}} \right) + f_{1,3} \right),
\end{align*}  
\begin{align*} 
\frac{\partial^{2} w_{1}}{\partial x_{1}^{2}} + \frac{\partial^{2} w_{1}}{\partial x_{2}^{2}} + \frac{\partial^{2} w_{1}} {\partial x_{3}^{2}}  &- \frac{1}{\left( \mathrm{c_{a}} + \mathrm{c_{d}} \right)} ( \frac{\partial w_{1}}{\partial t}  - \left( \mathrm{c_{0}} + \mathrm{c_{d}} - \mathrm{c_{a}} \right) \left( \frac{\partial^{2} w_{1}}{\partial x_{1}^{2}}  + \frac{\partial^{2} w_{2}}{\partial x_{1} \partial x_{2}} +  \frac{\partial^{2} w_{3}}{\partial x_{1} \partial x_{3}} \right) \\ &+ u_{1}\frac{\partial w_{1}}{\partial x_{1}} + u_{2}\frac{\partial w_{1}}{\partial x_{2}} + u_{3}\frac{\partial w_{1}}{\partial x_{3}} + 4 \mathrm{v_{r}}w_{1} )  = - \frac{1}{\left( \mathrm{c_{a}} + \mathrm{c_{d}} \right)} \left( 2 \mathrm{v_{r}} \left( \frac{\partial u_{3}}{\partial x_{2}} - \frac{\partial u_{2}}{\partial x_{3}} \right) + f_{2,1} \right) ,
\end{align*}  
\begin{align*}  
\frac{\partial^{2} w_{2}}{\partial x_{1}^{2}} + \frac{\partial^{2} w_{2}}{\partial x_{2}^{2}} + \frac{\partial^{2} w_{2}} {\partial x_{3}^{2}} &- \frac{1}{\left( \mathrm{c_{a}} + \mathrm{c_{d}} \right)} ( \frac{\partial w_{2}}{\partial t} - \left( \mathrm{c_{0}} + \mathrm{c_{d}} - \mathrm{c_{a}} \right) \left( \frac{\partial^{2} w_{1}}{\partial x_{2} \partial x_{1}} + \frac{\partial^{2} w_{2}}{\partial^{2} x_{2}} +  \frac{\partial^{2} w_{3}}{\partial x_{2} \partial x_{3}} \right) \\
&+ u_{1}\frac{\partial w_{2}}{\partial x_{1}} + u_{2}\frac{\partial w_{2}}{\partial x_{2}} + u_{3}\frac{\partial w_{2}}{\partial x_{3}} + 4 \mathrm{v_{r}}w_{2} )  = - \frac{1}{\left( \mathrm{c_{a}} + \mathrm{c_{d}} \right)} \left( 2 \mathrm{v_{r}} \left( \frac{\partial u_{1}}{\partial x_{3}} - \frac{\partial u_{3}}{\partial x_{1}} \right) + f_{2,2} \right),
\end{align*}  
\begin{align*} 
\frac{\partial^{2} w_{3}}{\partial x_{1}^{2}} + \frac{\partial^{2} w_{3}}{\partial x_{2}^{2}} + \frac{\partial^{2} w_{3}} {\partial x_{3}^{2}}  &- \frac{1}{\left( \mathrm{c_{a}} + \mathrm{c_{d}} \right)} ( \frac{\partial w_{3}}{\partial t} - \left( \mathrm{c_{0}} + \mathrm{c_{d}} - \mathrm{c_{a}} \right) \left( \frac{\partial^{2} w_{1}}{\partial x_{3} \partial x_{1}} + \frac{\partial^{2} w_{2}}{\partial x_{3} \partial x_{2}} +  \frac{\partial^{2} w_{3}}{\partial^{2} x_{3}} \right) \\
&+ u_{1}\frac{\partial w_{3}}{\partial x_{1}} + u_{2}\frac{\partial w_{3}}{\partial x_{2}} + u_{3}\frac{\partial w_{3}}{\partial x_{3}} + 4 \mathrm{v_{r}}w_{3} )  = - \frac{1}{\left( \mathrm{c_{a}} + \mathrm{c_{d}} \right)} \left( 2 \mathrm{v_{r}} \left( \frac{\partial u_{2}}{\partial x_{1}} - \frac{\partial u_{1}}{\partial x_{2}} \right) + f_{2,3} \right) ,
\end{align*} 
\begin{align*}
\frac{\partial^{2} p}{\partial x_{1}^{2}} + \frac{\partial^{2} p}{\partial x_{2}^{2}} + \frac{\partial^{2} p} {\partial x_{3}^{2}} &= - \frac{\partial}{\partial x_{1}} \left( u_{1}\frac{\partial u_{1}}{\partial x_{1}} + u_{2}\frac{\partial u_{1}}{\partial x_{2}} + u_{3}\frac{\partial u_{1}}{\partial x_{3}} + 2 \mathrm{v_{r}} \left( \frac{\partial w_{3}}{\partial x_{2}} - \frac{\partial w_{2}}{\partial x_{3}} \right) + f_{1,1} \right) \\
&- \frac{\partial}{\partial x_{2}} \left( u_{1}\frac{\partial u_{2}}{\partial x_{1}} + u_{2}\frac{\partial u_{2}}{\partial x_{2}} + u_{3}\frac{\partial u_{2}}{\partial x_{3}} + 2 \mathrm{v_{r}} \left( \frac{\partial w_{1}}{\partial x_{3}} - \frac{\partial w_{3}}{\partial x_{1}} \right) + f_{1,2} \right) \\
&- \frac{\partial}{\partial x_{3}} \left( u_{1}\frac{\partial u_{3}}{\partial x_{1}} + u_{2}\frac{\partial u_{3}}{\partial x_{2}} + u_{3}\frac{\partial u_{3}}{\partial x_{3}} + 2 \mathrm{v_{r}} \left( \frac{\partial w_{2}}{\partial x_{1}} - \frac{\partial w_{1}}{\partial x_{2}} \right) + f_{1,3} \right). 
\end{align*}
Then, the correctional functionals for above equations in $x_{1}-$direction are  
\begin{align*}
u_{1,n+1}(x_{1},x_{2},x_{3},t) &= u_{1,n}(x_{1},x_{2},x_{3},t) + \int_{\Omega_{x_{1}}} \lambda_{1}(X_{1}) \bigg[ \frac{\partial^{2} u_{1,n}(X_{1},x_{2},x_{3},t)}{\partial x_{1}^{2}} + \frac{\partial^{2} \tilde{u}_{1,n}(X_{1},x_{2},x_{3},t)}{\partial x_{2}^{2}} \\
&+ \frac{\partial^{2} \tilde{u}_{1,n}(X_{1},x_{2},x_{3},t)} {\partial x_{3}^{2}} - \frac{1}{\left( \mathrm{v} + \mathrm{v_{r}} \right)} \bigg( \frac{\partial \tilde{u}_{1,n}(X_{1},x_{2},x_{3},t)}{\partial t} + \tilde{u}_{1,n}(X_{1},x_{2},x_{3},t) \frac{\partial \tilde{u}_{1,n} (X_{1},x_{2},x_{3},t)}{\partial X_{1}} \\
&+ \tilde{u}_{2,n}(X_{1},x_{2},x_{3},t) \frac{\partial \tilde{u}_{1,n}(X_{1},x_{2},x_{3},t)}{\partial x_{2}} + \tilde{u}_{3,n}(X_{1},x_{2},x_{3},t) \frac{\partial \tilde{u}_{1,n}(X_{1},x_{2},x_{3},t)}{\partial x_{3}} \\
&+ \frac{\partial \tilde{p}_{n}(X_{1},x_{2},x_{3},t)}{\partial X_{1}} \bigg) + \frac{1}{\bigg( \mathrm{v} + \mathrm{v_{r}} \bigg)} \bigg( 2 \mathrm{v_{r}} \bigg( \frac{\partial \tilde{w}_{3,n}(X_{1},x_{2},x_{3},t)}{\partial x_{2}} - \frac{\partial \tilde{w}_{2,n}(X_{1},x_{2},x_{3},t)}{\partial x_{3}} \bigg) \\
&+ f_{1,1}(X_{1},x_{2},x_{3},t) \bigg) \bigg] d X_{1}, 
\end{align*}
\begin{align*} 
u_{2,n+1}(x_{1},x_{2},x_{3},t) &= u_{2,n}(x_{1},x_{2},x_{3},t) + \int_{\Omega_{x_{1}}} \lambda_{2}(X_{1}) \bigg[ \frac{\partial^{2} u_{2,n}(X_{1},x_{2},x_{3},t)}{\partial X_{1}^{2}} + \frac{\partial^{2} \tilde{u}_{2,n}(X_{1},x_{2},x_{3},t)}{\partial x_{2}^{2}} \\
&+ \frac{\partial^{2} \tilde{u}_{2,n}(X_{1},x_{2},x_{3},t)} {\partial x_{3}^{2}} - \frac{1}{\left( \mathrm{v} + \mathrm{v_{r}} \right)} \bigg( \frac{\partial \tilde{u}_{2,n}(X_{1},x_{2},x_{3},t)}{\partial t} + \tilde{u}_{1,n}(X_{1},x_{2},x_{3},t) \frac{\partial \tilde{u}_{2,n}(X_{1},x_{2},x_{3},t)}{\partial X_{1}} \\
&+ \tilde{u}_{2,n}(X_{1},x_{2},x_{3},t) \frac{\partial \tilde{u}_{2,n}(X_{1},x_{2},x_{3},t)}{\partial x_{2}} + \tilde{u}_{3,n}(X_{1},x_{2},x_{3},t) \frac{\partial \tilde{u}_{2,n}(X_{1},x_{2},x_{3},t)}{\partial x_{3}} \\
&+ \frac{\partial \tilde{p}_{n}(X_{1},x_{2},x_{3},t)}{\partial x_{2}} + \frac{1}{\bigg( \mathrm{v} + \mathrm{v_{r}} \bigg)} \bigg( 2 \mathrm{v_{r}} \bigg( \frac{\partial \tilde{w}_{1,n}(X_{1},x_{2},x_{3},t)}{\partial x_{3}} - \frac{\partial \tilde{w}_{3,n}(X_{1},x_{2},x_{3},t)}{\partial x_{1}} \bigg) \\
&+ f_{1,2}(X_{1},x_{2},x_{3},t) \bigg) \bigg] d X_{1}, 
\end{align*}
\begin{align*} 
u_{3,n+1}(x_{1},x_{2},x_{3},t) &= u_{3,n}(x_{1},x_{2},x_{3},t) + \int_{\Omega_{x_{1}}} \lambda_{3}(X_{1}) \bigg[ \frac{\partial^{2} u_{1,n}(X_{1},x_{2},x_{3},t)}{\partial X_{1}^{2}} + \frac{\partial^{2} \tilde{u}_{1,n}(X_{1},x_{2},x_{3},t)}{\partial x_{2}^{2}} \\
&+ \frac{\partial^{2} \tilde{u}_{3,n}(X_{1},x_{2},x_{3},t)} {\partial x_{3}^{2}} - \frac{1}{\left( \mathrm{v} + \mathrm{v_{r}} \right)} \bigg( \frac{\partial \tilde{u}_{3,n}(X_{1},x_{2},x_{3},t)}{\partial t} + \tilde{u}_{1,n}(X_{1},x_{2},x_{3},t) \frac{\partial \tilde{u}_{3,n}(X_{1},x_{2},x_{3},t)}{\partial X_{1}} \\
&+ \tilde{u}_{2,n}(X_{1},x_{2},x_{3},t) \frac{\partial \tilde{u}_{3,n}(X_{1},x_{2},x_{3},t)}{\partial x_{2}} + \tilde{u}_{3,n}(X_{1},x_{2},x_{3},t) \frac{\partial \tilde{u}_{3,n}(X_{1},x_{2},x_{3},t)}{\partial x_{3}} \\
&+ \frac{\partial \tilde{p}_{n}(X_{1},x_{2},x_{3},t)}{\partial x_{3}} \bigg) + \frac{1}{\left( \mathrm{v} + \mathrm{v_{r}} \right)} \bigg( 2 \mathrm{v_{r}} \bigg( \frac{\partial \tilde{w}_{2,n}(X_{1},x_{2},x_{3},t)}{\partial X_{1}} - \frac{\partial \tilde{w}_{1,n}(X_{1},x_{2},x_{3},t)}{\partial x_{2}} \bigg) \\
&+ f_{1,3}(X_{1},x_{2},x_{3},t) \bigg) \bigg]  d X_{1},
\end{align*}
\begin{align*}
w_{1,n+1}(x_{1},x_{2},x_{3},t) &= w_{1,n}(x_{1},x_{2},x_{3},t) + \int_{\Omega_{x_{1}}} \lambda_{4}(X_{1}) \bigg[ \frac{\partial^{2} w_{1,n}(X_{1},x_{2},x_{3},t)}{\partial X_{1}^{2}} + \frac{\partial^{2} \tilde{w}_{1,n}(X_{1},x_{2},x_{3},t)}{\partial x_{2}^{2}} \\
& + \frac{\partial^{2} \tilde{w}_{1,n}(X_{1},x_{2},x_{3},t)} {\partial x_{3}^{2}} - \frac{1}{ \left( \mathrm{c_{a}} + \mathrm{c_{d}} \right)} \bigg( \frac{\partial \tilde{w}_{1,n}(X_{1},x_{2},x_{3},t)}{\partial t} \\
&- \left( \mathrm{c_{0}} + \mathrm{c_{d}} - \mathrm{c_{a}} \right) \bigg( \frac{\partial^{2} \tilde{w}_{1,n}(X_{1},x_{2},x_{3},t)}{\partial X_{1}^{2}} + \frac{\partial^{2} \tilde{w}_{2,n}(X_{1},x_{2},x_{3},t)}{\partial X_{1} \partial x_{2}} +  \frac{\partial^{2} \tilde{w}_{3,n}(X_{1},x_{2},x_{3},t)}{\partial X_{1} \partial x_{3}} \bigg) \\
&+ \tilde{u}_{1,n}(X_{1},x_{2},x_{3},t) \frac{\partial \tilde{w}_{1,n}(X_{1},x_{2},x_{3},t)}{\partial X_{1}} + \tilde{u}_{2,n}(X_{1},x_{2},x_{3},t) \frac{\partial \tilde{w}_{1,n}(X_{1},x_{2},x_{3},t)}{\partial x_{2}} \\
&+ \tilde{u}_{3,n}(X_{1},x_{2},x_{3},t)\frac{\partial \tilde{w}_{1,n}(X_{1},x_{2},x_{3},t)}{\partial x_{3}} + 4 \mathrm{v_{r}} \tilde{w}_{1,n}(X_{1},x_{2},x_{3},t) \bigg)  \\
&+ \frac{1}{\left( \mathrm{c_{a}} + \mathrm{c_{d}} \right)} \left( 2 \mathrm{v_{r}} \left( \frac{\partial \tilde{u}_{3,n}(X_{1},x_{2},x_{3},t)}{\partial x_{2}} - \frac{\partial \tilde{u}_{2,n}(X_{1},x_{2},x_{3},t)}{\partial x_{3}} \right) + f_{2,1}(X_{1},x_{2},x_{3},t) \right) \bigg] d X_{1}, 
\end{align*}  
\begin{align*}  
w_{2,n+1}(x_{1},x_{2},x_{3},t) &= w_{2,n}(x_{1},x_{2},x_{3},t) + \int_{0}^{t} \lambda_{5}(X_{1}) \bigg[ \frac{\partial^{2} w_{2,n}(X_{1},x_{2},x_{3},t)}{\partial X_{1}^{2}} + \frac{\partial^{2} \tilde{w}_{2,n}(X_{1},x_{2},x_{3},t)}{\partial x_{2}^{2}} \\
&+ \frac{\partial^{2} \tilde{w}_{2,n}(X_{1},x_{2},x_{3},t)} {\partial x_{3}^{2}} - \frac{1}{\left( \mathrm{c_{a}} + \mathrm{c_{d}} \right)} \bigg( \frac{\partial \tilde{w}_{2,n}(X_{1},x_{2},x_{3},t)}{\partial t} \\
&- \left( \mathrm{c_{0}} + \mathrm{c_{d}} - \mathrm{c_{a}} \right) \bigg( \frac{\partial^{2} \tilde{w}_{1,n}(X_{1},x_{2},x_{3},t)}{\partial x_{2} \partial X_{1}} + \frac{\partial^{2} \tilde{w}_{2,n}(X_{1},x_{2},x_{3},t)}{\partial^{2} x_{2}} +  \frac{\partial^{2} \tilde{w}_{3,n}(X_{1},x_{2},x_{3},t)}{\partial x_{2} \partial x_{3}} \bigg) \\
&+ \tilde{u}_{1,n}(X_{1},x_{2},x_{3},t) \frac{\partial \tilde{w}_{2,n}(X_{1},x_{2},x_{3},t)}{\partial X_{1}} + \tilde{u}_{2,n}(X_{1},x_{2},x_{3},t) \frac{\partial \tilde{w}_{2,n}(X_{1},x_{2},x_{3},t)}{\partial x_{2}} \\
&+ \tilde{u}_{3,n}(X_{1},x_{2},x_{3},t) \frac{\partial \tilde{w}_{2,n}(X_{1},x_{2},x_{3},t)}{\partial x_{3}} + 4 \mathrm{v_{r}} \tilde{w}_{2,n}(X_{1},x_{2},x_{3},t) ) \\
&+ \frac{1}{\left( \mathrm{c_{a}} + \mathrm{c_{d}} \right)} \bigg( 2 \mathrm{v_{r}} \bigg( \frac{\partial \tilde{u}_{1,n}(X_{1},x_{2},x_{3},t)}{\partial x_{3}} - \frac{\partial \tilde{u}_{3,n}(X_{1},x_{2},x_{3},t)}{\partial X_{1}} \bigg) + f_{2,2}(X_{1},x_{2},x_{3},t) \bigg) \bigg] d X_{1} ,
\end{align*}  
\begin{align*} 
w_{3,n+1}(x_{1},x_{2},x_{3},t) &= w_{3,n}(x_{1},x_{2},x_{3},t) + \int_{0}^{t} \lambda_{6}(X_{1}) \bigg[ \frac{\partial^{2} w_{3,n}(X_{1},x_{2},x_{3},t)}{\partial X_{1}^{2}} + \frac{\partial^{2} \tilde{w}_{3,n}(X_{1},x_{2},x_{3},t)}{\partial x_{2}^{2}} \\
&+ \frac{\partial^{2} \tilde{w}_{3,n}(X_{1},x_{2},x_{3},t)} {\partial x_{3}^{2}} - \frac{1}{\left( \mathrm{c_{a}} + \mathrm{c_{d}} \right)} \bigg( \frac{\partial \tilde{w}_{3,n}(X_{1},x_{2},x_{3},t)}{\partial t} \\
&- \left( \mathrm{c_{0}} + \mathrm{c_{d}} - \mathrm{c_{a}} \right) \bigg( \frac{\partial^{2} \tilde{w}_{1,n}(X_{1},x_{2},x_{3},t)}{\partial x_{3} \partial X_{1}} + \frac{\partial^{2} \tilde{w}_{2,n}(X_{1},x_{2},x_{3},t)}{\partial x_{3} \partial x_{2}} +  \frac{\partial^{2} \tilde{w}_{3,n}(X_{1},x_{2},x_{3},t)}{\partial^{2} x_{3}} \bigg) \\
&+ \tilde{u}_{1,n}(X_{1},x_{2},x_{3},t) \frac{\partial \tilde{w}_{3,n}(X_{1},x_{2},x_{3},t)}{\partial X_{1}} + \tilde{u}_{2,n}(X_{1},x_{2},x_{3},t) \frac{\partial \tilde{w}_{3,n}(X_{1},x_{2},x_{3},t)}{\partial x_{2}} 
\\
&+ \tilde{u}_{3,n}(X_{1},x_{2},x_{3},t) \frac{\partial \tilde{w}_{3,n}(X_{1},x_{2},x_{3},t)}{\partial x_{3}} + 4 \mathrm{v_{r}} \tilde{w}_{3,n}(X_{1},x_{2},x_{3},t) \bigg) \\
&+ \frac{1}{\left( \mathrm{c_{a}} + \mathrm{c_{d}} \right)} \left( 2 \mathrm{v_{r}} \left( \frac{\partial \tilde{u}_{2,n}(X_{1},x_{2},x_{3},t)}{\partial X_{1}} - \frac{\partial \tilde{u}_{1,n}(X_{1},x_{2},x_{3},t)}{\partial x_{2}} \right) + f_{2,3}(X_{1},x_{2},x_{3},t) \right) \bigg] d X_{1} ,
\end{align*} 
\begin{align*}
p_{n+1}(x_{1},x_{2},x_{3},t) &= p_{n}(x_{1},x_{2},x_{3},t) + \int_{\Omega_{x_{1}}} \mu(X_{1}) \bigg[ \frac{\partial^{2} p_{n}(X_{1},x_{2},x_{3},t)}{\partial X_{1}^{2}} + \frac{\partial^{2} \tilde{p}_{n}(X_{1},x_{2},x_{3},t)}{\partial x_{2}^{2}} + \frac{\partial^{2} \tilde{p}_{n}(X_{1},x_{2},x_{3},t)} {\partial x_{3}^{2}} \\
&+ \frac{\partial}{\partial X_{1}} \bigg( \tilde{u}_{1,n}(X_{1},x_{2},x_{3},t) \frac{\partial \tilde{u}_{1,n}(X_{1},x_{2},x_{3},t)}{\partial X_{1}} + \tilde{u}_{2,n}(X_{1},x_{2},x_{3},t) \frac{\partial \tilde{u}_{1,n}(X_{1},x_{2},x_{3},t)}{\partial x_{2}} \\
&+ \tilde{u}_{3,n}(X_{1},x_{2},x_{3},t) \frac{\partial \tilde{u}_{1,n}(X_{1},x_{2},x_{3},t)}{\partial x_{3}} - 2 \mathrm{v_{r}} \left( \frac{\partial \tilde{w}_{3,n}(X_{1},x_{2},x_{3},\tau)}{\partial x_{2}} - \frac{\partial \tilde{w}_{2,n}(X_{1},x_{2},x_{3},\tau)}{\partial x_{3}} \right) \\
&- f_{1,1}(X_{1},x_{2},x_{3},t) \bigg) + \frac{\partial}{\partial x_{2}} \bigg( \tilde{u}_{1,n}(X_{1},x_{2},x_{3},t) \frac{\partial \tilde{u}_{2,n}(X_{1},x_{2},x_{3},t)}{\partial X_{1}} \\
&+ \tilde{u}_{2,n}(X_{1},x_{2},x_{3},t) \frac{\partial \tilde{u}_{2,n}(X_{1},x_{2},x_{3},t)}{\partial x_{2}} + \tilde{u}_{3,n}(X_{1},x_{2},x_{3},t) \frac{\partial \tilde{u}_{2,n}(X_{1},x_{2},x_{3},t)}{\partial x_{3}} \\
&- 2 \mathrm{v_{r}} \left( \frac{\partial \tilde{w}_{1,n}(X_{1},x_{2},x_{3},\tau)}{\partial x_{3}} - \frac{\partial \tilde{w}_{3,n}(X_{1},x_{2},x_{3},\tau)}{\partial x_{1}} \right) - f_{1,2}(X_{1},x_{2},x_{3},t) \bigg) \\
&+ \frac{\partial}{\partial x_{3}} \bigg( \tilde{u}_{1,n}(X_{1},x_{2},x_{3},t) \frac{\partial \tilde{u}_{3,n}(X_{1},x_{2},x_{3},t)}{\partial X_{1}} + \tilde{u}_{2,n}(X_{1},x_{2},x_{3},t) \frac{\partial \tilde{u}_{3,n}(X_{1},x_{2},x_{3},t)}{\partial x_{2}} \\
&+ \tilde{u}_{3,n}(X_{1},x_{2},x_{3},t) \frac{\partial \tilde{u}_{3,n}(X_{1},x_{2},x_{3},t)}{\partial x_{3}} - 2 \mathrm{v_{r}} \left( \frac{\partial \tilde{w}_{2,n}(\mathbf{x},\tau)}{\partial x_{1}} - \frac{\partial \tilde{w}_{1,n}(\mathbf{x},\tau)}{\partial x_{2}} \right) \\
&- f_{1,3}(X_{1},x_{2},x_{3},t) \bigg) \bigg] d X_{1}.
\end{align*}
Making the above correctional functionals stationary,
\begin{align*} 
\delta u_{1,n+1}(x_{1},x_{2},x_{3},t) &= \delta u_{1,n}(x_{1},x_{2},x_{3},t) + \left. \frac{\partial \delta u_{1,n}(X_{1},x_{2},x_{3},t)}{\partial X_{1}} \lambda_{1}(X_{1}) \right|_{X_{1} = x_{1}} - \left. \delta u_{1,n}(X_{1},x_{2},x_{3},t) \frac{d \lambda_{1}(X_{1})}{d X_{1}}  \right|_{X_{1} = x_{1}} \\
&+ \int_{\Omega_{x_{1}}} \delta u_{1,n}(X_{1},x_{2},x_{3},t) \frac{d^{2} \lambda_{1}(X_{1})}{d X^{2}_{1}} dX_{1} = 0,
\end{align*}  
\begin{align*} 
\delta u_{2,n+1}(x_{1},x_{2},x_{3},t) &= \delta u_{2,n}(x_{1},x_{2},x_{3},t) + \left. \frac{\partial \delta u_{2,n}(X_{1},x_{2},x_{3},t)}{\partial X_{1}} \lambda_{2}(X_{1}) \right|_{X_{1} = x_{1}} - \left. \delta u_{2,n}(X_{1},x_{2},x_{3},t) \frac{d \lambda_{2}(X_{1})}{d X_{1}}  \right|_{X_{1} = x_{1}} \\
&+ \int_{\Omega_{x_{1}}} \delta u_{2,n}(X_{1},x_{2},x_{3},t) \frac{d^{2} \lambda_{2}(X_{1})}{d X^{2}_{1}} dX_{1} = 0,
\end{align*}  
\begin{align*} 
\delta u_{3,n+1}(x_{1},x_{2},x_{3},t) &= \delta u_{3,n}(x_{1},x_{2},x_{3},t) + \left. \frac{\partial \delta u_{3,n}(X_{1},x_{2},x_{3},t)}{\partial X_{1}} \lambda_{3}(X_{1}) \right|_{X_{1} = x_{1}} - \left. \delta u_{3,n}(X_{1},x_{2},x_{3},t) \frac{d \lambda_{3}(X_{1})}{d X_{1}}  \right|_{X_{1} = x_{1}} \\
&+ \int_{\Omega_{x_{1}}} \delta u_{3,n}(X_{1},x_{2},x_{3},t) \frac{d^{2} \lambda_{3}(X_{1})}{d X^{2}_{1}} dX_{1} = 0, 
\end{align*} 
\begin{align*}
\delta w_{1,n+1}(x_{1},x_{2},x_{3},t) &= \delta w_{1,n}(x_{1},x_{2},x_{3},t) + \left. \frac{\partial \delta w_{1,n}(X_{1},x_{2},x_{3},t)}{\partial X_{1}} \lambda_{4}(X_{1}) \right|_{X_{1} = x_{1}} - \left. \delta w_{1,n}(X_{1},x_{2},x_{3},t) \frac{d \lambda_{4}(X_{1})}{d X_{1}}  \right|_{X_{1} = x_{1}} \\
&+ \int_{\Omega_{x_{1}}} \delta w_{1,n}(X_{1},x_{2},x_{3},t) \frac{d^{2} \lambda_{4}(X_{1})}{d X^{2}_{1}} dX_{1} = 0, 
\end{align*}  
\begin{align*} 
\delta w_{2,n+1}(x_{1},x_{2},x_{3},t) &= \delta w_{2,n}(x_{1},x_{2},x_{3},t) + \left. \frac{\partial \delta w_{2,n}(X_{1},x_{2},x_{3},t)}{\partial X_{1}} \lambda_{5}(X_{1}) \right|_{X_{1} = x_{1}} - \left. \delta w_{2,n}(X_{1},x_{2},x_{3},t) \frac{d \lambda_{5}(X_{1})}{d X_{1}}  \right|_{X_{1} = x_{1}} \\
&+ \int_{\Omega_{x_{1}}} \delta w_{2,n}(X_{1},x_{2},x_{3},t) \frac{d^{2} \lambda_{5}(X_{1})}{d X^{2}_{1}} dX_{1} = 0, 
\end{align*}  
\begin{align*} 
\delta w_{3,n+1}(x_{1},x_{2},x_{3},t) &= \delta w_{3,n}(x_{1},x_{2},x_{3},t) + \left. \frac{\partial \delta w_{3,n}(X_{1},x_{2},x_{3},t)}{\partial X_{1}} \lambda_{6}(X_{1}) \right|_{X_{1} = x_{1}} - \left. \delta w_{3,n}(X_{1},x_{2},x_{3},t) \frac{d \lambda_{6}(X_{1})}{d X_{1}}  \right|_{X_{1} = x_{1}} \\
&+ \int_{\Omega_{x_{1}}} \delta w_{3,n}(X_{1},x_{2},x_{3},t) \frac{d^{2} \lambda_{6}(X_{1})}{d X^{2}_{1}} dX_{1} = 0,  
\end{align*}  
\begin{align*}
\delta p_{n+1}(x_{1},x_{2},x_{3},t) &= \delta p_{n}(x_{1},x_{2},x_{3},t) + \left. \frac{\partial \delta p_{n}(X_{1},x_{2},x_{3},t)}{\partial X_{1}} \mu(X_{1}) \right|_{X_{1} = x_{1}} - \left. \delta p_{n}(X_{1},x_{2},x_{3},t) \frac{d \mu(X_{1})}{d X_{1}}  \right|_{X_{1} = x_{1}} \\
&+ \int_{\Omega_{x_{1}}} \delta p_{n}(X_{1},x_{2},x_{3},t) \frac{d^{2} \mu(X_{1})}{d X^{2}_{1}} dX_{1} = 0
\end{align*}
yields the following stationary conditions: for $i = 1,\cdots, 6$, that
\begin{align*}
\delta u_{i,n}: \frac{d^{2} \lambda_{i}}{d X^{2}_{1}} &= 0  & \ \ \ \delta p_{n}: \frac{d^{2} \mu}{d X^{2}_{1}} = 0 \\
\frac{\partial \delta u_{i,n}}{\partial X_{1}}: \left. \mu(X_{1}) \right|_{X_{1} = x_{1}} &= 1  \ \ \ \ \ \ \ \ \ \ \ \ \ \ \ \ \ \mathrm{and}& \ \ \ \frac{\partial \delta p_{n}}{\partial X_{1}}: \left. \mu(X_{1}) \right|_{X_{1} = x_{1}} = 1 \\
\delta u_{i,n}: \left. \frac{d \lambda_{i} (X_{1})}{d X_{1}} \right|_{X_{1} = x_{1}} &= 0  & \ \ \ \delta p_{n}: \left. \frac{d \mu (X_{1})}{d X_{1}} \right|_{X_{1} = x_{1}} = 0.
\end{align*}
Then, for $i = 1,\cdots, 6$, the Lagrange multipliers are $\lambda_{i} = 1$ and $\mu = 1$; and the desired iterative scheme is
\begin{align*}
u_{1,n+1}(x_{1},x_{2},x_{3},t) &= u_{1,n}(x_{1},x_{2},x_{3},t) + \int_{\Omega_{x_{1}}} \bigg[ \frac{\partial^{2} u_{1,n}(X_{1},x_{2},x_{3},t)}{\partial x_{1}^{2}} + \frac{\partial^{2} \tilde{u}_{1,n}(X_{1},x_{2},x_{3},t)}{\partial x_{2}^{2}} \\
&+ \frac{\partial^{2} \tilde{u}_{1,n}(X_{1},x_{2},x_{3},t)} {\partial x_{3}^{2}} - \frac{1}{\left( \mathrm{v} + \mathrm{v_{r}} \right)} \bigg( \frac{\partial \tilde{u}_{1,n}(X_{1},x_{2},x_{3},t)}{\partial t} + \tilde{u}_{1,n}(X_{1},x_{2},x_{3},t) \frac{\partial \tilde{u}_{1,n} (X_{1},x_{2},x_{3},t)}{\partial X_{1}} \\
&+ \tilde{u}_{2,n}(X_{1},x_{2},x_{3},t) \frac{\partial \tilde{u}_{1,n}(X_{1},x_{2},x_{3},t)}{\partial x_{2}} + \tilde{u}_{3,n}(X_{1},x_{2},x_{3},t) \frac{\partial \tilde{u}_{1,n}(X_{1},x_{2},x_{3},t)}{\partial x_{3}} \\
&+ \frac{\partial \tilde{p}_{n}(X_{1},x_{2},x_{3},t)}{\partial X_{1}} \bigg) + \frac{1}{\bigg( \mathrm{v} + \mathrm{v_{r}} \bigg)} \bigg( 2 \mathrm{v_{r}} \bigg( \frac{\partial \tilde{w}_{3,n}(X_{1},x_{2},x_{3},t)}{\partial x_{2}} - \frac{\partial \tilde{w}_{2,n}(X_{1},x_{2},x_{3},t)}{\partial x_{3}} \bigg) \\
&+ f_{1,1}(X_{1},x_{2},x_{3},t) \bigg) \bigg] d X_{1}, 
\end{align*}  
\begin{align*} 
u_{2,n+1}(x_{1},x_{2},x_{3},t) &= u_{2,n}(x_{1},x_{2},x_{3},t) + \int_{\Omega_{x_{1}}} \bigg[ \frac{\partial^{2} u_{2,n}(X_{1},x_{2},x_{3},t)}{\partial X_{1}^{2}} + \frac{\partial^{2} \tilde{u}_{2,n}(X_{1},x_{2},x_{3},t)}{\partial x_{2}^{2}} \\
&+ \frac{\partial^{2} \tilde{u}_{2,n}(X_{1},x_{2},x_{3},t)} {\partial x_{3}^{2}} - \frac{1}{\left( \mathrm{v} + \mathrm{v_{r}} \right)} \bigg( \frac{\partial \tilde{u}_{2,n}(X_{1},x_{2},x_{3},t)}{\partial t} + \tilde{u}_{1,n}(X_{1},x_{2},x_{3},t) \frac{\partial \tilde{u}_{2,n}(X_{1},x_{2},x_{3},t)}{\partial X_{1}} \\
&+ \tilde{u}_{2,n}(X_{1},x_{2},x_{3},t) \frac{\partial \tilde{u}_{2,n}(X_{1},x_{2},x_{3},t)}{\partial x_{2}} + \tilde{u}_{3,n}(X_{1},x_{2},x_{3},t) \frac{\partial \tilde{u}_{2,n}(X_{1},x_{2},x_{3},t)}{\partial x_{3}} \\
&+ \frac{\partial \tilde{p}_{n}(X_{1},x_{2},x_{3},t)}{\partial x_{2}} + \frac{1}{\bigg( \mathrm{v} + \mathrm{v_{r}} \bigg)} \bigg( 2 \mathrm{v_{r}} \bigg( \frac{\partial \tilde{w}_{1,n}(X_{1},x_{2},x_{3},t)}{\partial x_{3}} - \frac{\partial \tilde{w}_{3,n}(X_{1},x_{2},x_{3},t)}{\partial x_{1}} \bigg) \\
&+ f_{1,2}(X_{1},x_{2},x_{3},t) \bigg) \bigg] d X_{1}, 
\end{align*}  
\begin{align*} 
u_{3,n+1}(x_{1},x_{2},x_{3},t) &= u_{3,n}(x_{1},x_{2},x_{3},t) + \int_{\Omega_{x_{1}}} \bigg[ \frac{\partial^{2} u_{1,n}(X_{1},x_{2},x_{3},t)}{\partial X_{1}^{2}} + \frac{\partial^{2} \tilde{u}_{1,n}(X_{1},x_{2},x_{3},t)}{\partial x_{2}^{2}} \\
&+ \frac{\partial^{2} \tilde{u}_{3,n}(X_{1},x_{2},x_{3},t)} {\partial x_{3}^{2}} - \frac{1}{\left( \mathrm{v} + \mathrm{v_{r}} \right)} \bigg( \frac{\partial \tilde{u}_{3,n}(X_{1},x_{2},x_{3},t)}{\partial t} + \tilde{u}_{1,n}(X_{1},x_{2},x_{3},t) \frac{\partial \tilde{u}_{3,n}(X_{1},x_{2},x_{3},t)}{\partial X_{1}} \\
&+ \tilde{u}_{2,n}(X_{1},x_{2},x_{3},t) \frac{\partial \tilde{u}_{3,n}(X_{1},x_{2},x_{3},t)}{\partial x_{2}} + \tilde{u}_{3,n}(X_{1},x_{2},x_{3},t) \frac{\partial \tilde{u}_{3,n}(X_{1},x_{2},x_{3},t)}{\partial x_{3}} \\
&+ \frac{\partial \tilde{p}_{n}(X_{1},x_{2},x_{3},t)}{\partial x_{3}} \bigg) + \frac{1}{\left( \mathrm{v} + \mathrm{v_{r}} \right)} \bigg( 2 \mathrm{v_{r}} \bigg( \frac{\partial \tilde{w}_{2,n}(X_{1},x_{2},x_{3},t)}{\partial X_{1}} - \frac{\partial \tilde{w}_{1,n}(X_{1},x_{2},x_{3},t)}{\partial x_{2}} \bigg) \\
&+ f_{1,3}(X_{1},x_{2},x_{3},t) \bigg) \bigg] d X_{1},
\end{align*} 
\begin{align*}
w_{1,n+1}(x_{1},x_{2},x_{3},t) &= w_{1,n}(x_{1},x_{2},x_{3},t) + \int_{\Omega_{x_{1}}} \bigg[ \frac{\partial^{2} w_{1,n}(X_{1},x_{2},x_{3},t)}{\partial X_{1}^{2}} + \frac{\partial^{2} \tilde{w}_{1,n}(X_{1},x_{2},x_{3},t)}{\partial x_{2}^{2}} \\
& + \frac{\partial^{2} \tilde{w}_{1,n}(X_{1},x_{2},x_{3},t)} {\partial x_{3}^{2}} - \frac{1}{ \left( \mathrm{c_{a}} + \mathrm{c_{d}} \right)} \bigg( \frac{\partial \tilde{w}_{1,n}(X_{1},x_{2},x_{3},t)}{\partial t} \\
&- \left( \mathrm{c_{0}} + \mathrm{c_{d}} - \mathrm{c_{a}} \right) \bigg( \frac{\partial^{2} \tilde{w}_{1,n}(X_{1},x_{2},x_{3},t)}{\partial X_{1}^{2}} + \frac{\partial^{2} \tilde{w}_{2,n}(X_{1},x_{2},x_{3},t)}{\partial X_{1} \partial x_{2}} +  \frac{\partial^{2} \tilde{w}_{3,n}(X_{1},x_{2},x_{3},t)}{\partial X_{1} \partial x_{3}} \bigg) \\
&+ \tilde{u}_{1,n}(X_{1},x_{2},x_{3},t) \frac{\partial \tilde{w}_{1,n}(X_{1},x_{2},x_{3},t)}{\partial X_{1}} + \tilde{u}_{2,n}(X_{1},x_{2},x_{3},t) \frac{\partial \tilde{w}_{1,n}(X_{1},x_{2},x_{3},t)}{\partial x_{2}} \\
&+ \tilde{u}_{3,n}(X_{1},x_{2},x_{3},t)\frac{\partial \tilde{w}_{1,n}(X_{1},x_{2},x_{3},t)}{\partial x_{3}} + 4 \mathrm{v_{r}} \tilde{w}_{1,n}(X_{1},x_{2},x_{3},t) \bigg)  \\
&+ \frac{1}{\left( \mathrm{c_{a}} + \mathrm{c_{d}} \right)} \left( 2 \mathrm{v_{r}} \left( \frac{\partial \tilde{u}_{3,n}(X_{1},x_{2},x_{3},t)}{\partial x_{2}} - \frac{\partial \tilde{u}_{2,n}(X_{1},x_{2},x_{3},t)}{\partial x_{3}} \right) + f_{2,1}(X_{1},x_{2},x_{3},t) \right) \bigg] d X_{1}, 
\end{align*}  
\begin{align*} 
w_{2,n+1}(x_{1},x_{2},x_{3},t) &= w_{2,n}(x_{1},x_{2},x_{3},t) + \int_{0}^{t} \bigg[ \frac{\partial^{2} w_{2,n}(X_{1},x_{2},x_{3},t)}{\partial X_{1}^{2}} + \frac{\partial^{2} \tilde{w}_{2,n}(X_{1},x_{2},x_{3},t)}{\partial x_{2}^{2}} \\
&+ \frac{\partial^{2} \tilde{w}_{2,n}(X_{1},x_{2},x_{3},t)} {\partial x_{3}^{2}} - \frac{1}{\left( \mathrm{c_{a}} + \mathrm{c_{d}} \right)} \bigg( \frac{\partial \tilde{w}_{2,n}(X_{1},x_{2},x_{3},t)}{\partial t} \\
&- \left( \mathrm{c_{0}} + \mathrm{c_{d}} - \mathrm{c_{a}} \right) \bigg( \frac{\partial^{2} \tilde{w}_{1,n}(X_{1},x_{2},x_{3},t)}{\partial x_{2} \partial X_{1}} + \frac{\partial^{2} \tilde{w}_{2,n}(X_{1},x_{2},x_{3},t)}{\partial^{2} x_{2}} +  \frac{\partial^{2} \tilde{w}_{3,n}(X_{1},x_{2},x_{3},t)}{\partial x_{2} \partial x_{3}} \bigg) \\
&+ \tilde{u}_{1,n}(X_{1},x_{2},x_{3},t) \frac{\partial \tilde{w}_{2,n}(X_{1},x_{2},x_{3},t)}{\partial X_{1}} + \tilde{u}_{2,n}(X_{1},x_{2},x_{3},t) \frac{\partial \tilde{w}_{2,n}(X_{1},x_{2},x_{3},t)}{\partial x_{2}} \\
&+ \tilde{u}_{3,n}(X_{1},x_{2},x_{3},t) \frac{\partial \tilde{w}_{2,n}(X_{1},x_{2},x_{3},t)}{\partial x_{3}} + 4 \mathrm{v_{r}} \tilde{w}_{2,n}(X_{1},x_{2},x_{3},t) ) \\
&+ \frac{1}{\left( \mathrm{c_{a}} + \mathrm{c_{d}} \right)} \bigg( 2 \mathrm{v_{r}} \bigg( \frac{\partial \tilde{u}_{1,n}(X_{1},x_{2},x_{3},t)}{\partial x_{3}} - \frac{\partial \tilde{u}_{3,n}(X_{1},x_{2},x_{3},t)}{\partial X_{1}} \bigg) + f_{2,2}(X_{1},x_{2},x_{3},t) \bigg) \bigg] d X_{1} ,
\end{align*}  
\begin{align*} 
w_{3,n+1}(x_{1},x_{2},x_{3},t) &= w_{3,n}(x_{1},x_{2},x_{3},t) + \int_{0}^{t} \bigg[ \frac{\partial^{2} w_{3,n}(X_{1},x_{2},x_{3},t)}{\partial X_{1}^{2}} + \frac{\partial^{2} \tilde{w}_{3,n}(X_{1},x_{2},x_{3},t)}{\partial x_{2}^{2}} \\
&+ \frac{\partial^{2} \tilde{w}_{3,n}(X_{1},x_{2},x_{3},t)} {\partial x_{3}^{2}} - \frac{1}{\left( \mathrm{c_{a}} + \mathrm{c_{d}} \right)} \bigg( \frac{\partial \tilde{w}_{3,n}(X_{1},x_{2},x_{3},t)}{\partial t} \\
&- \left( \mathrm{c_{0}} + \mathrm{c_{d}} - \mathrm{c_{a}} \right) \bigg( \frac{\partial^{2} \tilde{w}_{1,n}(X_{1},x_{2},x_{3},t)}{\partial x_{3} \partial X_{1}} + \frac{\partial^{2} \tilde{w}_{2,n}(X_{1},x_{2},x_{3},t)}{\partial x_{3} \partial x_{2}} +  \frac{\partial^{2} \tilde{w}_{3,n}(X_{1},x_{2},x_{3},t)}{\partial^{2} x_{3}} \bigg) \\
&+ \tilde{u}_{1,n}(X_{1},x_{2},x_{3},t) \frac{\partial \tilde{w}_{3,n}(X_{1},x_{2},x_{3},t)}{\partial X_{1}} + \tilde{u}_{2,n}(X_{1},x_{2},x_{3},t) \frac{\partial \tilde{w}_{3,n}(X_{1},x_{2},x_{3},t)}{\partial x_{2}} 
\\
&+ \tilde{u}_{3,n}(X_{1},x_{2},x_{3},t) \frac{\partial \tilde{w}_{3,n}(X_{1},x_{2},x_{3},t)}{\partial x_{3}} + 4 \mathrm{v_{r}} \tilde{w}_{3,n}(X_{1},x_{2},x_{3},t) \bigg) \\
&+ \frac{1}{\left( \mathrm{c_{a}} + \mathrm{c_{d}} \right)} \left( 2 \mathrm{v_{r}} \left( \frac{\partial \tilde{u}_{2,n}(X_{1},x_{2},x_{3},t)}{\partial X_{1}} - \frac{\partial \tilde{u}_{1,n}(X_{1},x_{2},x_{3},t)}{\partial x_{2}} \right) + f_{2,3}(X_{1},x_{2},x_{3},t) \right) \bigg] d X_{1} ,
\end{align*} 
\begin{align*}
p_{n+1}(x_{1},x_{2},x_{3},t) &= p_{n}(x_{1},x_{2},x_{3},t) + \int_{\Omega_{x_{1}}} \bigg[ \frac{\partial^{2} p_{n}(X_{1},x_{2},x_{3},t)}{\partial X_{1}^{2}} + \frac{\partial^{2} \tilde{p}_{n}(X_{1},x_{2},x_{3},t)}{\partial x_{2}^{2}} + \frac{\partial^{2} \tilde{p}_{n}(X_{1},x_{2},x_{3},t)} {\partial x_{3}^{2}} \\
&+ \frac{\partial}{\partial X_{1}} \bigg( \tilde{u}_{1,n}(X_{1},x_{2},x_{3},t) \frac{\partial \tilde{u}_{1,n}(X_{1},x_{2},x_{3},t)}{\partial X_{1}} + \tilde{u}_{2,n}(X_{1},x_{2},x_{3},t) \frac{\partial \tilde{u}_{1,n}(X_{1},x_{2},x_{3},t)}{\partial x_{2}} \\
&+ \tilde{u}_{3,n}(X_{1},x_{2},x_{3},t) \frac{\partial \tilde{u}_{1,n}(X_{1},x_{2},x_{3},t)}{\partial x_{3}} - 2 \mathrm{v_{r}} \left( \frac{\partial \tilde{w}_{3,n}(X_{1},x_{2},x_{3},\tau)}{\partial x_{2}} - \frac{\partial \tilde{w}_{2,n}(X_{1},x_{2},x_{3},\tau)}{\partial x_{3}} \right) \\
&- f_{1,1}(X_{1},x_{2},x_{3},t) \bigg) + \frac{\partial}{\partial x_{2}} \bigg( \tilde{u}_{1,n}(X_{1},x_{2},x_{3},t) \frac{\partial \tilde{u}_{2,n}(X_{1},x_{2},x_{3},t)}{\partial X_{1}} \\
&+ \tilde{u}_{2,n}(X_{1},x_{2},x_{3},t) \frac{\partial \tilde{u}_{2,n}(X_{1},x_{2},x_{3},t)}{\partial x_{2}} + \tilde{u}_{3,n}(X_{1},x_{2},x_{3},t) \frac{\partial \tilde{u}_{2,n}(X_{1},x_{2},x_{3},t)}{\partial x_{3}} \\
&- 2 \mathrm{v_{r}} \left( \frac{\partial \tilde{w}_{1,n}(X_{1},x_{2},x_{3},\tau)}{\partial x_{3}} - \frac{\partial \tilde{w}_{3,n}(X_{1},x_{2},x_{3},\tau)}{\partial x_{1}} \right) - f_{1,2}(X_{1},x_{2},x_{3},t) \bigg) \\
&+ \frac{\partial}{\partial x_{3}} \bigg( \tilde{u}_{1,n}(X_{1},x_{2},x_{3},t) \frac{\partial \tilde{u}_{3,n}(X_{1},x_{2},x_{3},t)}{\partial X_{1}} + \tilde{u}_{2,n}(X_{1},x_{2},x_{3},t) \frac{\partial \tilde{u}_{3,n}(X_{1},x_{2},x_{3},t)}{\partial x_{2}} \\
&+ \tilde{u}_{3,n}(X_{1},x_{2},x_{3},t) \frac{\partial \tilde{u}_{3,n}(X_{1},x_{2},x_{3},t)}{\partial x_{3}} - 2 \mathrm{v_{r}} \left( \frac{\partial \tilde{w}_{2,n}(\mathbf{x},\tau)}{\partial x_{1}} - \frac{\partial \tilde{w}_{1,n}(\mathbf{x},\tau)}{\partial x_{2}} \right) \\
&- f_{1,3}(X_{1},x_{2},x_{3},t) \bigg) \bigg] d X_{1}.
\end{align*}
\section{Conclusion}
He's variational iteration method is generalized to create an iterative scheme that gurantee a unique solution for asymmetric systems of PDEs with respect to the time. Systems of incompressible fluid flow and incompressible micropolar fluid flow can be solved uniquely by applying the generalized He's variational iteration method to create a desired iterative scheme for them. 

\paragraph{Acknowledgement:} The author received no direct funding for this research.

\end{document}